\documentclass{amsart}

\newtheorem{thm}{Theorem}[section]
\newtheorem{theorem}{Theorem}[section]
\newtheorem{lemma}[thm]{Lemma}

\newtheorem{example}[thm]{Example}
\newtheorem{question}[thm]{Question}

\theoremstyle{definition}

\def\eqn#1$$#2$${\begin{equation}\label#1#2\end{equation}}

\def\th{{\leavevmode\setbox1=\hbox{t}
  \hbox to \wd1{t\kern-0.6ex{\char039}\hss}}}
\def\dh{{\leavevmode\setbox1=\hbox{d}
  \hbox to 1.05\wd1{d\kern-0.4ex{\char039}\hss}}}
\def\=#1{\if #1u{\accent23u}\else
\ifx #1d{\dh}\else \ifx #1t{\th}\else
 {\accent20 #1}\fi\fi\fi}
\def\'#1{\if #1i{\accent19\i}\else {\accent19 #1}\fi}

\def\diam{\operatorname{diam}}

\def\en{\mathbb N}
\def\er{\mathbb R}

\def\dist{\operatorname{dist}}

\def \reg {\partial _{\kern1pt\text{reg}}}

\newtoks\by
\newtoks\paper
\newtoks\book
\newtoks\jour
\newtoks\yr
\newtoks\pages
\newtoks\vol
\newtoks\publ
\newtoks\eds
\newtoks\proc
\newtoks\mathrev
\newtoks\web
\def\ota{{\hbox{???}}}
\def\cLear{\by=\ota\paper=\ota\book=\ota\jour=\ota\yr=\ota
\pages=\ota\vol=\ota\publ=\ota}
\def\endpaper{\the\by, \textit{\the\paper},
{\the\jour} \textbf{\the\vol} (\the\yr), \the\pages.\cLear}
\def\endbook{\the\by, \textit{\the\book}, \the\publ, \the\yr.\cLear}
\def\endprep{\the\by, \textit{\the\paper}, \the\jour.}
\def\endprepkma{\the\by, \textit{\the\paper}, \the\jour, available on http://adela.\-karlin.\-mff.\-cuni.\-cz/\~{ }rokyta/\-preprint/\-index.php.\cLear}
\def\endproc{\the\by, \textit{\the\paper}, in:\,\the\book,
\the\publ, \the\yr, \the\pages.\cLear}
\def\endper{\the\by, \textit{personal communication}.\cLear}

\begin{document}
\title{Baire-one mappings contained in a usco map}
\author{Ond\v{r}ej F.K. Kalenda}

\address{Charles University\\Faculty of Mathematics and Physics\\
Department of Mathematical Analysis\\
Sokolovsk\'a~83\\
186~75 Praha~8\\
Czech Republic}

\email{kalenda@karlin.mff.cuni.cz}
\thanks{The work is a part of the research project MSM 0021620839 financed by MSMT and partly supported by the research grant GA \v{C}R 201/06/0018 and
by Universidad Polit\'ecnica de Valencia.}



\begin{abstract} We investigate Baire--one functions whose graph is contained in a graph of usco mapping. We prove
in particular that such a function defined on a metric space with
values in $\mathbb{R}^d$ is the pointwise limit of a sequence of
continuous functions with graphs contained in the graph of a common
usco map.
\end{abstract}

\maketitle

\section{Introduction}

We study the following question:

\bigskip

{\it Let $X$ be a metric space, $Y$ a convex subset of a normed
linear space and $f:X\to Y$ a Baire-one function whose graph is
contained in the graph of a usco mapping. Is there a sequence of
continuous functions $f_n:X\to Y$ pointwise converging to $f$ such
that the graphs of all the $f_n$'s are contained in a usco map
$\varphi:X\to Y$?}

\bigskip

This question  appeared in the joint research of R.Anguelov and the
author (see \cite{AK}). Moreover, this question is also natural and
interesting in itself as there are several theorems on the existence
of Baire-one selections of multivalued (in particular usco) maps --
see e.g. \cite{HJT,JR,srivatsa}. In particular, every usco map from
a metric space into a normed linear space admits a Baire--one
selection; therefore Baire-one functions whose graph is contained in
the graph of a usco mapping are quite common.

We do not know the full answer to our question but we prove some
partial results. One of them is that the answer is positive if $Y$
is a closed convex subset of a finite-dimensional space. This is
used in \cite{AK} to show that if $X$ is a Baire metric space then
continuous functions from $X$ to $\er^d$ form a dense subset of the
convergence space of minimal usco maps.

Let us start by recalling and introducing some notions.

A nonempty-valued mapping $\varphi:X\to Y$ is called \textit{upper
semi-continuous compact valued} (shortly \textit{usco}) if
$\varphi(x)$ is a (nonempty) compact subset of $Y$ for each $x\in X$
and $\{x\in X:\varphi(x)\subset U\}$ is open in $X$ for each
$U\subset Y$ open.

A function $f:X\to Y$ is called \textit{Baire-one} if it is the
pointwise limit of a sequence of continuous functions.

We say that a family of functions (defined on $X$ with values in
$Y$) is \textit{usco-bounded} if there is a usco map $\varphi:X\to
Y$ whose graph, i.e. the set
$$\{(x,y)\in X\times Y: y\in \varphi(x)\}$$
contains the graphs of all the functions from the family. We will
use this terminology for single functions and for sequences of
functions.

\section{Basic facts and examples}

An important role is played by the following characterization of
usco maps and maps whose graph is contained in the graph of a usco
map.

\begin{lemma}\label{charusco}
Let $X$ and $Y$ be metric spaces and $\varphi:X\to Y$ a
nonempty-valued set-valued mapping.
\begin{itemize}
    \item[(i)] The mapping $\varphi$ is usco if and only if whenever $x_n$
    is a sequence in $X$ converging to some $x\in X$ and $y_n\in
    \varphi(x_n)$ for each $n\in\mathbb{N}$, there is a subsequence of
    $y_n$ converging to an element of $\varphi(x)$.
    \item[(ii)] There is a usco map $\psi:X\to Y$ with
    $\varphi\subset\psi$ (in the sense of inclusion of graphs)
    if and only if whenever $x_n$
    is a sequence in $X$ converging to some $x\in X$ and $y_n\in
    \varphi(x_n)$ for each $n\in\mathbb{N}$, there is a convergent subsequence of
    $y_n$.
\end{itemize}
\end{lemma}

\begin{proof}
The point (i) is an analogue of \cite[Lemma 3.1.1]{fabian}. The
assertion of the quoted lemma is the same -- only $X$ and $Y$ are
arbitrary topological spaces and nets are used instead of sequences.
As we are now dealing with metric spaces, sequences are enough and
the same proof works.

Let us show the point (ii). The `only if' part follows immediately
from (i). It remains to prove the `if' part. Let $\psi$ be the
multivalued mapping whose graph is the closure in $X\times Y$ of the
graph of $\varphi$. We will show that $\psi$ is usco using (i).
Denote by $d$ and $\rho$ the metrics of $X$ and $Y$, respectively.

Let $x_n$ be a sequence in $X$ converging to $x\in X$ and
$y_n\in\psi(x_n)$ for each $n\in\en$. Then each pair $(x_n,y_n)$
belongs to the graph of $\psi$. As the graph of $\varphi$ is dense
in the graph of $\psi$, there are pairs $(x'_n,y'_n)$ in the graph
of $\varphi$ such that $d(x'_n,x_n)<\frac1n$ and
$\rho(y'_n,y_n)<\frac1n$ for each $n\in\en$. Then $x'_n\to x$, and
hence there is a subsequence $y'_{n_k}$ converging to some $y\in Y$.
Then $y_{n_k}$ converge to $y$ as well. Hence $(x_{n_k},y_{n_k})$
converges to $(x,y)$. As the graph of $\psi$ is closed, we get
$y\in\psi(x)$. This completes the proof.
\end{proof}

An important subclass of Baire--one functions consists of so-called
simple functions. We recall the definition.

Let $X$ and $Y$ be  metric spaces and $f:X\to Y$ be a function. The
function $f$ is called \textit{simple} if there is a
$\sigma$-discrete partition of $X$ into $F_\sigma$ sets such that
$f$ is constant on each element of the partition.

Recall that a family of subsets of $X$ is \textit{discrete} if each
point of $X$ has a neighborhood meeting at most one element of the
family; and a family is \textit{$\sigma$-discrete} if it is a
countable union of discrete families.

It is easy to check that $f$ is simple if and only if there is a
partition $F_\gamma$, $\gamma\in\Gamma$, of $X$ such that
\begin{itemize}
    \item $f$ is constant on each $F_\gamma$.
    \item Each $F_\gamma$ can be expressed as an increasing union of
    closed sets $F_\gamma^n$ such that the family $F_\gamma^n$,
    $\gamma\in\Gamma$, is discrete for each $n\in\mathbb{N}$.
\end{itemize}

If $Y$ is a convex subset of a normed linear space (or more
generally if $Y$ is arcwise connected metric space), then any simple
function with values in $Y$ is Baire--one (see e.g. \cite[Lemma
2.13]{sel}).

\bigskip

Moreover, we need the following well-known approximation result.

\begin{lemma}\label{b1approx} Let $X$ and $Y$ be metric spaces. Then
any Baire--one function $f:X\to Y$ is the uniform limit of a
sequence of simple functions.
\end{lemma}

\begin{proof} It follows from \cite[Lemmata 1.1 and 1.4]{hansell} that $f$  has a $\sigma$-discrete function base consisting of closed sets
(i.e., a $\sigma$-discrete family $\mathcal{B}$ of closed subsets of
$X$ such that $f^{-1}(U)$ is the union of a subfamily of
$\mathcal{B}$ for each $U\subset Y$ open). The assertion then
follows by using Lemma 2.7 of \cite{sel}.
\end{proof}

Now we are going to give some examples of Baire--one mappings which
are usco-bounded and some which are not. The first one is trivial.
It follows from the fact that the closed bounded sets in $\er^n$ are
compact.

\begin{example} Any bounded (Baire-one) function $f:X\to \er^n$ is
usco-bounded. \end{example}

Also the second example is trivial:

\begin{example} (i) The function $f:\er\to\er$ defined by
$$f(x)=\begin{cases} \frac1x, & x\ne 0, \\ 0, & x=0.\end{cases}$$ is a
Baire-one function which is not usco-bounded.

(ii) The function $f:\er\to\er$ defined by $f(x)=x$ is usco-bounded
(it is continuous) although it is not bounded.
\end{example}

We continue by two more general statements.

\begin{example} Let $Y$ be an infinite dimensional normed space.
Then there is a bounded Baire-one function $f:\er\to Y$ which is not
usco-bounded.
\end{example}

\begin{proof} Let $y_n$ be a sequence in the unit ball of $Y$ which
has no convergent subsequence (as $Y$ is infinite-dimensional, the
unit ball is not compact and hence such a sequence exists). Define
the function $f$ by the formula
$$f(x)=\begin{cases} 0, & x\in (-\infty,0]\cup(1,+\infty), \\ y_n, &
x\in(\frac1{n+1},\frac1n],\,n\in\en.\end{cases}$$
Then $f$ is bounded, it is easily seen to be simple, and hence
Baire-one. However, it is not usco-bounded due to
Lemma~\ref{charusco}. Indeed, $\frac1n\to0$ and the sequence
$f(\frac1n)=y_n$ has no convergent subsequence.
\end{proof}

\begin{example}\label{protipriklad} Let $Y$ be a non-complete metric space. Then
there is a sequence of simple functions $f_n:\er\to Y$ which
uniformly converges to a simple function $f:\er\to Y$ such that each
$f_n$ is usco-bounded but $f$ is not.
\end{example}

\begin{proof}
Let $y_n$ be a Cauchy sequence in $Y$ which is not convergent.
Define the functions $f_n$ by the formula:
$$ f_n(x)=\begin{cases}0, & x\in (-\infty,0]\cup(1,+\infty), \\ y_n, &
x\in(0,\frac1n], \\ y_k, & x\in(\frac1{k+1},\frac1k],\,1\le
k<n.\end{cases}$$

Then each $f_n$ is easily seen to be a usco-bounded simple function.
Moreover, as the sequence $y_n$ is Cauchy, the sequence $f_n$
uniformly converges to the function $f$ defined by the formula
$$f(x)=\begin{cases} 0, & x\in (-\infty,0]\cup(1,+\infty), \\ y_n, &
x\in(\frac1{n+1},\frac1n],\,n\in\en.\end{cases}$$
This function is simple and is not usco-bounded. This can be proved
by  the argument used in the previous example (note that $y_n$ is a
non-converging Cauchy sequence and hence has no convergent
subsequence).
\end{proof}

\section{Main results}

In this section we show some partial positive answers to the
question from the introduction. The first result shows that the
answer is positive if $f$ is a simple function.

\begin{theorem}\label{simple}
Let $X$ be a metric space, $Y$ a convex subset of a normed linear
space and $f:X\to Y$ be a usco-bounded simple function.  Then there
is a usco-bounded sequence of continuous functions $f_n:X\to Y$
which pointwise converges to $f$.
\end{theorem}

\begin{proof}
We will imitate the proof of the fact that simple functions are
Baire--one in \cite[Lemmata 2.12 and 2.13]{sel}. Without loss of
generality suppose that $0\in Y$. Fix a partition $F_\gamma$,
$\gamma\in\Gamma$, of $X$ such that
\begin{itemize}
    \item $f$ is constant on each $F_\gamma$.
    \item Each $F_\gamma$ can be expressed as an increasing union of
    closed sets $F_\gamma^n$ such that the family $F_\gamma^n$,
    $\gamma\in\Gamma$, is discrete for each $n\in\mathbb{N}$.
\end{itemize}

For $n\in\mathbb{N}$ and $\gamma\in\Gamma$ we define the following
functions:
\begin{align*}
    d^n(x)&=\dist(x,\bigcup_{\gamma\in\Gamma} F^n_\gamma), \\
    d^n_\gamma(x)&=\dist(x,F^n_\gamma), \\
    e^n_\gamma(x)&=\dist(x,\bigcup_{\delta\in\Gamma\setminus\{\gamma\}}
    F^n_\delta).
\end{align*}
All these functions are continuous and hence the set
$$G^n_\gamma=\{ x\in X: d^n_\gamma(x)<\tfrac13 e^n_\gamma(x)\}$$
is open for each $n$ and $\gamma$. Moreover, as the family
$F_\gamma^n$, $\gamma\in\Gamma$, is discrete, we get
$F_\gamma^n\subset G^n_\gamma$. It is proved in \cite[pp.
34--35]{sel} that the family $G_\gamma^n$, $\gamma\in\Gamma$, is
discrete for each $n\in\mathbb{N}$, too.

Denote by $y_\gamma$ the value of $f$ on $F_\gamma$ and define the
functions $f_n$ as follows:
$$ f_n(x)=\begin{cases} 0, & x\in X\setminus\bigcup_{\gamma\in\Gamma} G^n_\gamma, \\
\left(1-\frac{4d^n_\gamma(x)}{d^n_\gamma(x)+e^n_\gamma(x)}\right)y_\gamma,
& x\in G^n_\gamma.\end{cases}$$ The function $f_n$ is continuous
(see \cite[pp. 35--36]{sel}) and $f_n(x)=y_\gamma$ for $x\in
F^n_\gamma$. Hence clearly $f_n$ pointwise converges to $f$.

It remains to show that the sequence $f_n$ is usco-bounded. To see
this we will use Lemma \ref{charusco}. Fix a sequence $x_n$
converging to $x\in X$ and a sequence $k_n$ of natural numbers. We
will be done if we show that the sequence $f_{k_n}(x_n)$ has a
converging subsequence.

Up to passing to a subsequence we can suppose that the sequence
$k_n$ is either constant or increasing. If $k_n=k$ for each $n$,
then $f_{k_n}(x_n)=f_k(x_n)$ converges to $f_k(x)$ due to the
continuity of $f_k$.

Hence suppose that $k_n$ is increasing. If $x_n\in
X\setminus\bigcup_{\gamma\in\Gamma} G^{k_n}_\gamma$ for infinitely
many $n$'s, then $f_{k_n}(x_n)=0$ for infinitely many $n$'s and
hence we have converging subsequence.

Thus suppose that for each $n\in\mathbb{N}$ there is some
$\gamma_n\in\Gamma$ with $x_n\in G^{k_n}_{\gamma_n}$. By the
definition of the set $G^{k_n}_{\gamma_n}$ there is some $z_n\in
F^{k_n}_{\gamma_n}$ with
$\rho(x_n,z_n)<\frac13e^{k_n}_{\gamma_n}(x_n)$. We have
$f_{k_n}(x_n)=c_n y_{\gamma_n}$ for some $c_n\in[0,1]$. Without loss
of generality we may suppose that the sequence $c_n$ converges to
some $c\in[0,1]$.

Let $\alpha\in\Gamma$ be such that $x\in F_\alpha$. If
$\gamma_n=\alpha$ for infinitely many $n$'s, then there is a
subsequence of $f_{k_n}(x_n)$ converging to $c y_\alpha$.

Finally, suppose that $\gamma_n\ne\alpha$ for each $n$. Then we have
$$
\rho(x_n,z_n)<\frac13e^{k_n}_{\gamma_n}(x_n) \le
\frac13e^{k_n}_{\gamma_n}(x)+\frac13\rho(x_n,x)
$$
for $n\in\en$. As there is some $n_0$ such that $x\in
F^{n_0}_\alpha$, $e^{k_n}_{\gamma_n}(x)=0$ for $n$ sufficiently
large. Therefore $\rho(x_n,z_n)$ converges to $0$ and thus $z_n$
converges to $x$. As $f$ is usco-bounded we get (by Lemma
\ref{charusco}) a converging subsequence of $f(z_n)=y_{\gamma_n}$.
Then $f_{k_n}(x_n)=c_n y_{\gamma_n}$ has converging subsequence as
well.
\end{proof}

The next result is an analogue of the standard fact that Baire--one
functions are preserved by the uniform limits. Note, that the
assumption that the range is complete is necessary, while in the
standard setting completeness is not needed.

\begin{theorem}\label{limita} Let $X$ be a metric space and $Y$ be a convex subset of a normed linear
space which is complete in the norm metric. Let $f_n:X\to Y$ be a
sequence of mappings which uniformly converges to a mapping $f:X\to
Y$. If each $f_n$ is the pointwise limit of a usco-bounded sequence
of continuous functions, then $f$ has the same property.

The completeness assumption on $Y$ cannot be omitted.
\end{theorem}

\begin{proof} We will imitate the proof of the fact that
Baire-one functions are preserved by uniform limits given in
\cite[Lemma 2.14]{sel}. As the quoted proof contains a large number
of misprints, we give a complete proof.

Denote by $Z$ the normed space which contains $Y$. Without loss of
generality suppose that $\|f_m(x)-f(x)\|<2^{-m}$ for each $x\in X$
and $m\in\mathbb{N}$. For each $m\in \mathbb{N}$ let
$(f_{m,n}:n\in\mathbb{N})$ be a usco-bounded sequence pointwise
converging to $f_m$.

We define continuous functions $g_{m,n}:X\to Y$ as follows:
\begin{gather*}
    g_{1,n}= f_{1,n} \quad \mbox{for }n\in\mathbb{N}, \\
    g_{m+1,n}(x)=\begin{cases} f_{m+1,n}(x), &
    \|f_{m+1,n}(x)-g_{m,n}(x)\|\le 2^{-m+1}, \\
    g_{m,n}(x)+
    2^{-m+1}\frac{f_{m+1,n}(x)-g_{m,n}(x)}{\|f_{m+1,n}(x)-g_{m,n}(x)\|},
    & \|f_{m+1,n}(x)-g_{m,n}(x)\| > 2^{-m+1}. \end{cases}
\end{gather*}

We have
\begin{equation*}
    \forall x\in X\,\forall m\in\mathbb{N}\,\exists
    n_0\in\mathbb{N}\,\forall n\ge n_0: g_{m,n}(x)=f_{m,n}(x).
\end{equation*}
Indeed, let $x\in X$ be arbitrary. As $g_{1,n}= f_{1,n}$ for all
$n$, the assertion is true for $m=1$. Suppose it is true for some
$m$. As
$$\|f_m(x)-f_{m+1}(x)\|\le\|f_m(x)-f(x)\|+\|f(x)-f_{m+1}(x)\|<
2^{-m}+2^{-m-1}<2^{-m+1},$$ we have
$\|f_{m,n}(x)-f_{m+1,n}(x)\|<2^{-m+1}$ for $n$ large enough. Now
using the induction hypothesis and the definition of $g_{m+1,n}(x)$,
we get $g_{m+1,n}(x)=f_{m+1,n}(x)$ for $n$ large enough.

Set $h_n=g_{n,n}$ for $n\in\mathbb{N}$. Then the sequence $h_n$
pointwise converges to $f$ and, moreover, is usco-bounded.

Let us show the first assertion. Let $x\in X$ and $\varepsilon>0$ be
arbitrary. Choose $m\in\mathbb{N}$ such that $2^{-m+4}<\varepsilon$.
Fix $n_0\in\mathbb{N}$ such that $g_{m,n}(x)=f_{m,n}(x)$ and
$\|f_{m,n}(x)-g_m(x)\|<\frac\varepsilon2$ for $n\ge n_0$. Then for
$n\ge n_0$ we have
\begin{align*}
    \|h_n(x)-f(x)\|&\le
    \|g_{n,n}(x)-g_{m,n}(x)\|+\|g_{m,n}(x)-f_m(x)\|+\|f_m(x)-f(x)\|
    \\ & \le 2^{-m+1}+\dots+2^{-n+1} +\|f_{m,n}(x)-f_m(x)\|+2^{-m}
    \\ & < 2^{-m+3}+\frac\varepsilon2<\varepsilon.
\end{align*}
This shows that $h_n(x)$ converges to $f(x)$.

Now we are going to prove that the sequence $h_n$ is usco-bounded.
So take an arbitrary sequence $x_n\in X$ converging to some $x\in X$
and a sequence $k_n$ of natural numbers. We need to show that the
sequence $h_{k_n}(x_n)$ has a converging subsequence.

If the sequence $k_n$ has a constant subsequence, we are done (as
$x_n\to x$ and each $f_k$ is continuous). Otherwise we can without
loss of generality suppose that the sequence $k_n$ is increasing.

Now, as $(f_{m,n}:n\in\mathbb{N})$ is usco-bounded for each
$m\in\mathbb{N}$, the sequence $f_{m,k_n}(x_n)$ has a converging
subsequence for each $m\in\mathbb{N}$. Thus, we can suppose without
loss of generality that, for each $m\in\mathbb{N}$ the sequence
$f_{m,k_n}(x_n)$ converges to some $y_m\in Y$.

Further, for each $n\in\mathbb{N}$ we have
$$ h_{k_n}(x_n)=f_{1,k_n}(x_n)+\sum_{j=1}^{k_n-1}
c_j^n\frac{f_{j+1,k_n}(x_n)-g_{j,k_n}(x_n)}{\|f_{j+1,k_n}(x_n)-g_{j,k_n}(x_n)\|},
$$
where $c_j^n\in[0,2^{-j+1}]$ for $n\in\mathbb{N}$ and
$j=1,\dots,k_n-1$. (If $f_{j+1,k_n}(x_n)=g_{j,k_n}(x_n)$, we set
$c_j^n=0$ and suppose the fraction equals to some unit vector.)

We can
consider sequences $(c_j^n : j=1,\dots,k_n-1)$ as elements of the
set
$$C=\{(t_j : j\in\mathbb{N}) : \forall j\in\mathbb{N}:
t_j\in[0,2^{-j+1}]\}.$$ The embedding is done by completing the
finite sequence by zeros since the $k_n$-th place. The set $C$ is a
compact subset of the Banach space $\ell_1$, hence we can suppose
without loss of generality that the sequences $(c_j^n)$ converge
(for $n\to\infty$) in the $\ell_1$-norm to a sequence
$(c_j:j\in\mathbb{N})\in C$.

Observe that that for each $j\in\mathbb{N}$ the sequence
$g_{j,k_n}(x_n)$ converges in $Y$. Indeed, for $j=1$ we have
$g_{1,k_n}(x_n)=f_{1,k_n}(x_n)$ which converges to $y_1$. Suppose
now that $j\in\mathbb{N}$ is such that $g_{j,k_n}(x_n)$ converges in
$Y$. Then
$$g_{j+1,k_n}(x_n)=g_{j,k_n}(x_n)+c_j^n\frac{f_{j+1,k_n}(x_n)-g_{j,k_n}(x_n)}{\|f_{j+1,k_n}(x_n)-g_{j,k_n}(x_n)\|},
$$
and hence it converges in $Z$ due to the assumption that
$f_{j+1,k_n}(x_n)$ converges to $y_{j+1}$ and that $c_j^n$ converges
to $c_j$. (If $f_{j+1,k_n}(x_n)=g_{j,k_n}(x_n)$ only for a finite
number of $n$'s, the conclusion is clear. If the equality holds for
infinitely many $n$'s then by the above convention $c_j^n=0$ for
infinitely many $n$'s, and hence $c_j=0$. So $c_j^n\to0$ and hence
the limit of $g_{j+1,k_n}(x_n)$ is the same as that of
$g_{j,k_n}(x_n)$.) Moreover, as $Y$ is complete and the values of
each $g_{m,n}$ are in $Y$, we get that the limit belongs to $Y$.

Finally, we have that for each $j\in\mathbb{N}$ the sequence
$c_j^n\frac{f_{j+1,k_n}(x_n)-g_{j,k_n}(x_n)}{\|f_{j+1,k_n}(x_n)-g_{j,k_n}(x_n)\|}$
converges for $n\to\infty$ to some $z_j\in Z$. As $\|z_j\|\le
2^{-j+1}$, the sequence $(y_1+\sum_{j=1}^{N-1} z_j:N\in\en)$ is
Cauchy. Moreover, $y_1+\sum_{j=1}^{N-1} z_j\in Y$ for each
$N\in\en$, as it is equal to
$$\lim_{n\to\infty} \left(f_{1,k_n}(x_n)+\sum_{j=1}^{N-1}
c_j^n\frac{f_{j+1,k_n}(x_n)-g_{j,k_n}(x_n)}{\|f_{j+1,k_n}(x_n)-g_{j,k_n}(x_n)\|}\right)
=\lim_{n\to\infty} g_{N,k_n}(x_n),$$
which belongs to $Y$ as $Y$ is complete and $g_{N,k_n}(x_n)\in Y$
for each $N$ and $n$. Hence the above sequence converges in $Y$ to
$y_1+\sum_{j=1}^\infty z_j$. Now it is clear that $h_{k_n}(x_n)$
converges to $y_1+\sum_{j=1}^\infty z_j$, which completes the proof
of the positive part.

If $Y$ is not complete, Example~\ref{protipriklad} shows that the
assertion is not true. Indeed, there are usco-bounded simple
functions $f_n:\er\to Y$ uniformly converging to a function $f$
which is not usco-bounded. Each $f_n$ is the pointwise limit of a
usco-bounded sequence of continuous function by
Theorem~\ref{simple}. On the other hand, $f$ cannot be expressed as
such a limit as in such a case it would be usco-bounded.
\end{proof}

Finally, we give the result for the case when $Y$ is
finite-dimensional.

\begin{theorem}\label{Rn} Let $X$ be a metric space, $Y$ a closed convex subset of $\er^d$ and $f:X\to Y$ be a
usco-bounded Baire-one function. Then $f$ is the pointwise limit of
a usco-bounded sequence of continuous functions.
\end{theorem}

In the proof we will need two simple lemmata.

\begin{lemma}\label{modif} Let $X$ be a metric space, $f$ and $g$ two functions
defined on $X$ with values in $\mathbb{R}^d$ such that $f-g$ is
bounded. If $f$ is usco-bounded, then so is $g$.
\end{lemma}

\begin{proof} Let $x_n$ be a sequence in $X$ converging to some
$x\in X$. As $f$ is usco-bounded, we can suppose that the sequence
$f(x_n)$ converges. Further, the sequence $g(x_n)-f(x_n)$ is
bounded, and hence has a convergent subsequence. It follows that
$g(x_n)$ has a convergent subsequence. This shows that $g$ is
usco-bounded.
\end{proof}

\begin{lemma}\label{approx} Let $X$ be a metric space, $Y$ be a closed convex subset of $\er^d$ and $f:X\to
Y$ be a usco-bounded Baire-one function. Then $f$ is a uniform limit
of a sequence of usco-bounded simple functions.
\end{lemma}

\begin{proof} Fix $\varepsilon>0$. By Lemma \ref{b1approx} there is a simple
function $g$ with $\|f(x)-g(x)\|<\varepsilon$ for all $x\in X$. By
Lemma~\ref{modif} the function $g$ is usco-bounded. This completes
the proof.
\end{proof}

{\it Proof of Theorem~\ref{Rn}. } Let $f:X\to Y$ be a usco-bounded
Baire--one function. By Lemma~\ref{approx} it is the uniform limit
of a usco-bounded sequence of simple functions. Now the result
follows from Theorems~\ref{simple} and~\ref{limita}. \qed

\section{Final remarks and open questions}

Of course, the main problem is whether the answer to the question
from the introduction is positive in general. However, let us
formulate some more questions.

\begin{question} Let $X$ and $Y$ be metric spaces and $f:X\to Y$ a
usco-bounded Baire--one function. Is there a sequence of
usco-bounded simple functions uniformly converging to $f$?
\end{question}

The positive answer to this question would imply (using Theorem
\ref{limita}) the positive answer to our main problem under the
assumption that $Y$ is complete. Due to Example \ref{protipriklad}
it would not help to solve the general case.

By Lemma \ref{approx} the aswer to the above question is positive if
$Y$ is a closed convex subset of $\er^d$. Moreover, the answer is
positive if $f$ is continuous (and $X$, $Y$ are general metric
spaces). Let us sketch the proof of this (although it yields nothing
new with respect to our main problem).

Let $f:X\to Y$ be continuous and $\varepsilon>0$. Then for each
$x\in X$ there is an open neighborhood $U_x$ of $x$ with $\diam
f(U_x)<\varepsilon$. Let $\mathcal{V}$ be a locally finite open
refinement of the open cover $\{U_x:x\in X\}$ of $X$. Let
$\{V_\alpha:\alpha<\kappa\}$ be an enumeration of $\mathcal V$ by
ordinal numbers. Set
$$W_\alpha=V_\alpha\setminus\bigcup_{\beta<\alpha}V_\beta$$
for each $\alpha<\kappa$. If $W_\alpha\ne\emptyset$ choose
$y_\alpha\in f(W_\alpha)$. Define the function $g:X\to Y$ by setting
$g(x)= y_\alpha$ for $x\in W_\alpha$. It is clear that the distance
of $f(x)$ and $g(x)$ is less than $\varepsilon$ for each $x\in X$.
Further, it is easy to see that $g$ is a simple function. Finally,
as the partition $\{W_\alpha:\alpha<\kappa\}$ is locally finite, the
function $g$ is easily seen to be usco-bounded (using Lemma
\ref{charusco}).

Another question concerns possible modification of Theorem
\ref{limita}.

\begin{question} Let $X$ be a metric space and  $Y$ a convex subset
of a normed linear space. Suppose that the sequence of functions
$f_n:X\to Y$ is usco-bounded and uniformly converges to a function
$f$. Suppose, moreover, that each $f_n$ is the pointwise limit of a
usco-bounded sequence of continuous functions. Is the same true for
$f$? \end{question}

If $Y$ is complete, the answer is positive (even without the
assumption of the usco-boundedness of the sequence) by Theorem
\ref{limita}. Note that the sequence from Example \ref{protipriklad}
is not usco-bounded.

The positive answer to this question would help to solve the problem
if the following strengthening of the first question has positive
answer.

\begin{question} Let $X$ and $Y$ be metric spaces and $f:X\to Y$ a
usco-bounded Baire--one function. Is there a usco-bounded sequence
of simple functions uniformly converging to $f$?
\end{question}

Note, that the answer is positive if $Y$ is a closed convex subset
of $\er^d$. This follows easily from the proofs of Lemmata
\ref{modif} and \ref{approx}. (Note that if $f$ is usco-bounded and
$Y$ is a closed subset of $\er^n$, then the graph of the set-valued
map $x\mapsto B(f(x),\varepsilon)$ is contained in the graph of a
usco map -- we can use Lemma \ref{charusco}.)


\begin{thebibliography}{00}

\bibitem{AK}
\by{R.Anguelov and O.Kalenda} \paper{The convergence space of
minimal usco mappings} \jour{preprint}
\endprep

\bibitem{fabian}
\by{M. Fabian} \book{G\^ateaux differentiability of convex functions
and topology: weak Asplund spaces} \publ{Wiley-Interscience, New
York} \yr{1997}
\endbook

\bibitem{hansell}
\by{R.W. Hansell} \paper{First class functions with values in
nonseparable spaces} \book{Constantin Carath\'eodory: an
international tribute, Vol. I, II} \publ{World Sci. Publishing,
Teaneck, NJ} \yr{1991} \pages{461--475}
\endproc

\bibitem{HJT}
\by{R.W. Hansell, J.E. Jayne and M.Talagrand} \paper{First class
selectors for weakly upper semi-continuous multi-valued maps in
Banach spaces} \jour{J. Reine Angew. Math.} \vol{361} \yr{1985}
\pages{201--220}
\endpaper

\bibitem{JR}
\by{J.E. Jayne and  C.A. Rogers} \paper{Borel selectors for upper
semi-continuous set-valued mappings} \jour{Acta Math.} \vol{155}
\yr{1985} \pages{41--79}
\endpaper

\bibitem{sel}
\by{J.E. Jayne and  C.A. Rogers} \book{Selectors} \publ{Princeton
University Press} \yr{2002}
\endbook

\bibitem{srivatsa}
\by{V.V.  Srivatsa} \paper{Baire class 1 selectors for upper
semicontinuous set-valued maps} \jour{Trans. Amer. Math. Soc.}
\vol{337} \yr{1993} \pages{no. 2, 609--624}
\endpaper

\end{thebibliography}
\end{document}